\documentclass[preprint,12pt]{elsarticle}
\usepackage{amssymb}

\journal{Frontiers in Computational Neuroscience }

\begin{document}

\begin{frontmatter}

\title{Energy demands of diverse spiking cells from the neocortex, hippocampus and thalamus.}

\author{A. Moujahid, A. D'Anjou and M. Gra\~na}
\address{Computational Intelligence Group, Department of Computer Science and Artificial Intelligence, University of the Basque Country UPV/EHU, 20018 San Sebastian, Spain}

\date{September 15, 2013}

\begin{abstract}
It has long been known that neurons in the brain are not physiologically homogeneous. In response to current stimulus, they can fire several distinct patterns of action potentials that are associated with different physiological classes ranging from regular-spiking cells, fast-spiking cells, intrinsically bursting cells, and low-threshold cells. In this work we show that the high degree of variability in firing characteristics of action potentials among these cells is accompanied with a significant variability in the energy demands required to restore the concentration gradients after an action potential.
The values of the metabolic energy were calculated for a wide range of cell temperatures and stimulus intensities following two different approaches. The first one is based on the amount of Na$^+$ load crossing the membrane during a single action potential, while the second one focuses on the electrochemical energy functions deduced from the dynamics of the computational neuron models. The results show that the thalamocortical relay neuron is the most energy-efficient cell consuming between 7 to 18 nJ/cm$^2$ for each spike generated, while both the regular and fast spiking cells from somatosensory cortex and the intrinsically-bursting cell from a cat visual cortex are the least energy-efficient, and can consume up to 100 nJ/cm$^2$ per spike. The lowest values of these energy demands were achieved at higher temperatures and high external stimuli.

\end{abstract}



\end{frontmatter}

\maketitle

\section{Introduction}
A huge number of studies have been devoted to categorize neurons in the brain according to various criteria such as morphology, physiology, biochemical properties and synaptic characteristics \cite{Cauli1997, Peters1984, Toledo-Rodriguez2003, White1989}. However, relatively little has been said about their classification in terms of energy demands. The relationship between the action potentials characteristics and the energy cost required to generate them has been reported recently \cite{Sengupta2010,Carter2009,Crotty2006}. Nonetheless, as in many other works, the calculations of the metabolic energy involved in the generation of action potentials are based on the ion-counting methods that generally focus on the number of ATP molecules hydrolyzed by the sodium pump.

It is well known that Na$^+$ and K$^+$  are the most important ions with greater impact on the membrane potential dynamics and on its energetics \cite{Ames2000}. In fact, these two ions together with Ca$^{2+}$ account for the major ions whose movements consume ATP molecules. During the generation of action potentials, the electrochemical gradients  are partially altered and must be restored by Na$^+$/K$^+$ pump which mediates influx of 2 K$^+$  in exchange for 3 Na$^+$  per 1 ATP molecule consumed \cite{Kandel1991}. Estimations of the energy consumed by the neuron during its signaling activity are usually intended by extrapolating from a Hodgkin-Huxley type model the number of sodium ions required to depolarize the membrane in order to know the amount of ATP molecules that will be needed for the pump to reestablish its rest potential. This method, known as ion counting approach, corresponds to the first approach to account with the neuron energy  \cite{Laughlin1998,Alle2009,Attwell2001}. Taking into account that action potential propagation following Hodgkin/Huxley kinetics is based on equal Na$^+$ efflux and K$^+$ influx, this approach could lead to overestimate/underestimate values of energy \cite{Hertz2013}. For one hand, because of the asymmetry between Na$^+$,K$^+$-ATPase mediated Na$^+$ and K$^+$ fluxes. And, for the other hand, due to the overlap of inward and outward currents during an action potential generation. In fact, inward $Na^{+}$ and outward $K^{+}$ overlap during the action potential generation introducing an uncertainty in the calculation of sodium ions up to a factor of $4$ \cite{Attwell2001, Hodgkin1975, Lennie2003}.

Based on computational models that capture the intrinsic properties of action potentials from different mammalian neurons, this work provides a quantitative comparison of the amount of ionic currents crossing the membrane during an action potential and their corresponding energy demands focusing on a new approach. This approach is based on  the biophysical considerations about the nature of computational models used to account for neuronal spiking response, and gives an alternative method to deduce the electrochemical energy involved in the dynamics of the neuron model. Unlike the ion counting approach, this method requires no hypothesis about the stoichiometry of the ions or the extent of the overlapping between Na$^+$ and K$^+$ ions and, most importantly, offers the possibility to compute the energy consumption of coupled neurons \cite{Moujahid2011,Moujahid2010}, especially when considering complex networks where all the interactions should be considered. As commented in Ref. \cite{Buzsaki} it is not possible to estimate energetic costs in isolation. As far as possible the complexity of the interactions should be taken into account.

The analytical expression of the electrochemical energy involved in the dynamics of the classical Hodgkin-Huxley model consisting of sodium, potassium and leakage currents has been reported previously in \cite{Moujahid2011,Moujahid2012}. Here we provide energy functions of different spiking neurons where additional ionic currents are involved.  The intrinsic properties of the computational models, considered in this work, arise from voltage-dependent conductances described by Hodgkin-Huxley type kinetics. These models were adjusted in previous studies  \cite{Pospischil2008,Guo2008,Wang1996} to experimental data from different preparations to reproduce the firing characteristics of different neuron classes present in neocortex, thalamus and hippocampus. These neurons include regular spiking (RS) cells, fast spiking, intrinsically bursting (IB), hippocampal fast-spiking interneuron (HFI) and thalamocortical relay neuron (TCR). The RS and FS models can reproduce respectively firing properties of RS and FS cells as observed both in intracellular recordings of cells in ferret visual cortex and in rat somatosensory cortex. Meanwhile, IB neuron models reproduce the typical firings of intrinsically bursting cells in guinea-pig somatosensory cortex and in cat visual cortex. The other two neuron models capture the typical spiking as seen in rat hippocampal interneuron and a mouse thalamocortical relay neuron respectively.

The paper is organized as follows. First, we introduce the dynamics and electrochemical energy of the Hodgkin-Huxley like models used to account for the spiking response of different classes of neurons present in neocortex, thalamus and hippocampus. Then we describe the kinetic of these neurons according to the ionic currents involved in their dynamics. Results are reported and discussed in Section 3. Finally, Section 4 draws conclusions and open questions.

\section{Material and methods}
\subsection{Dynamics of the spiking neuron models}
All computational models were modeled using a Hodgkin-Huxley type of kinetic model \cite{Hodgkin1952}, and run using the Matlab simulation environment. For all these models, the dynamics governing the membrane potential obeys the following equation:

\begin{equation}
\begin{array}{ll}
C \, \dot{V}= \small -g_{l}(V-E_{l})-I_{Na}-I_K-I_M-I_L-I_T+I_{Stim},\\
\end{array}
\label{equ1}
\end{equation}
where $V$ is the membrane potential in mV, $g_{l}$ is the maximal conductance per unit area for the leakage channel, and $E_{l}$ is the corresponding reversal potential. $I_{Stim}$ is the total membrane current density in $\mu \textrm{A/cm}^2$. $C$ is the membrane capacitance in ($\mu$F/cm$^2$). $I_{Na}$, $I_K$ ($I_M$) and $I_L$ ($I_T$) are the sodium, potassium (slow potassium) and calcium (low-threshold calcium) currents respectively.

\begin{table}
\begin{tabular}{lll}

\hline
\hline
  & \small Activation &  \small  Inactivation  \\ \hline
 \small  Neocortical   &  &   \\
 \small cells  &  &   \\ \hline\\
$\scriptstyle I_{Na}=g_{Na}m^3h(V-E_{Na})$  & $ \scriptstyle\alpha_{m}(V)=\frac{-0.32(V-V_T-13)}{e^{-(V-V_T-13)/4}-1}$  & $ \scriptstyle\alpha_{h}(V)= \scriptstyle{0.128}e^{-(V-V_T-17)/18}$ \\
  & $ \scriptstyle\beta_{m}(V)=\frac{0.28(V-V_T-40)}{e^{(V-V_T-40)/5)}-1}$  & $ \scriptstyle\beta_{h}(V)=\frac{4}{e^{-(V-V_T-40)/5}+1}$ \\
$\scriptstyle I_{K}=g_Kn^4(V-E_K)$  & $ \scriptstyle\alpha_{n}(V)=\frac{-0.032(V-V_T-15)}{e^{-(V-V_T-15)/5}-1}$  &  \\
         & $ \scriptstyle\beta_{n}(V)=0.5 e^{-(V-V_T-10)/40}$  &  \\
$\scriptstyle I_M=g_Mp(V-E_K)$  & $\scriptstyle p_{\infty}(V)= \scriptstyle\frac{1}{e^{-(V+35)/10}+1}$  & $\scriptstyle\tau_p(V)= \scriptstyle\frac{\tau_{max}}{3.3e^{(V+35)/20}+e^{-(V+35)/20}}$ \\
$\scriptstyle I_{L}=g_{L}q^2r(V-E_{Ca})$  & $\scriptstyle\alpha_{q}(V)=\frac{0.055(-27-V)}{e^{(-27-V)/3.8}-1}$  &  $ \scriptstyle\alpha_{r}(V)= \scriptstyle{0.000457} e^{(-13-V)/50}$ \\
         & $ \scriptstyle\beta_{q}(V)= \scriptstyle{0.94} e^{(-75-V)/17}$  & $ \scriptstyle\beta_{r}(V)=\frac{0.0065}{e^{(-15-V)/28}+1}$ \\
\hline
 \small Thalamocortical  &  &   \\
 \small relay (TCR) cell  &  &   \\ \hline\\
$\scriptstyle I_{Na}=g_{Na}m_{\infty}^3h(V-E_{Na})$ & $\scriptstyle m_\infty(V)=\scriptstyle 1/(e^{-(V+37)/7}+1)$ & $\scriptstyle h_\infty(V)=1/(e^{(V+41)/4}+1)$\\
$\scriptstyle I_K=g_K(0.75(1-h))^4(V-E_K)$ &  & $\scriptstyle\tau_h(V)=1/(a_1(V)+b_1(V))$\\
&&$\scriptstyle a_1(V)=\scriptstyle 0.128 e^{(-(V+46)/18)}$ \\
&& $\scriptstyle b_1(V)=\scriptstyle 4/(1+e^{(-(V+23)/5)})$\\
$\scriptstyle I_T=g_Tp_{\infty}^2r(V-E_T)$&$\scriptstyle p_\infty(V)=\scriptstyle 1/(e^{-(V+60)/6.2}+1)$&$\scriptstyle r_\infty(V)=\scriptstyle 1/(e^{(V+84)/4}+1)$\\
&&$\scriptstyle \tau_r(V)=\scriptstyle 0.4(e^{-(V+25)/10.5}+28)$\\\\
\hline
 \small Hippocampal  &  &   \\
 \small interneuron (RHI)  &  &   \\ \hline\\
$\scriptstyle I_{Na}=g_{N{a}}m_{\infty}^{3}h(V-E_{N{a}})$  & $\scriptstyle m_{\infty}(V)=\scriptstyle \alpha_m(V)/(\alpha_m(V)+\beta_m(V))$  & \\
          & $\scriptstyle\alpha_m(V)=\scriptstyle \frac{-0.1(V+35)}{(e^{-0.1(V+35)}-1)}$  & $\scriptstyle\alpha_h(V)=\scriptstyle 0.07e^{-(V+58)/20}$\\
          &$\scriptstyle\beta_m(V)=\scriptstyle 4 e^{-(V+60)/18}$& $\scriptstyle\beta_h(V)=\scriptstyle1/(e^{-0.1(V+28)}+1)$\\
$\scriptstyle I_K=g_Kn^4(V-E_K)$  & $\scriptstyle\alpha_n(V)=\scriptstyle \frac{-0.01(V+34)}{(e^{-0.1(V+34)}-1)}$  &  \\
         & $\scriptstyle\beta_n(V)=\scriptstyle 0.125e^{-(V+44)/80}$  &  \\ \\

\hline\hline

\end{tabular}
\caption{The activation and inactivation functions describing ion currents. }
\label{table1}
\end{table}

\begin{table}
{\scriptsize
\begin{center}
\begin{tabular}{l||ccc|cc|ccc||c||c}
\hline\hline
                                & \multicolumn{8}{c||}{}  	          	     	  	 	   	   	        &    &  \\
                                & \multicolumn{8}{c||}{\bf Neocortex}  	          	     	& {\bf Thalamus}  & {\bf Hippocampus}\\
                                & \multicolumn{8}{c||}{}  	          	     	  	 	   	   	        &    &  \\\hline
                                & \multicolumn{3}{c|}{\bf RS cells} & \multicolumn{2}{c|}{\bf FS cells} & \multicolumn{3}{c||}{\bf IB cells} & {\bf TCR}  & {\bf RHI}  \\
                                &  	    &       & 	    & 	    & 	 	 &  	 &  	 &       &      & \\ \hline
{\bf Maximal}	                & Cell 1& Cell 2& Cell 3& Cell 4& Cell 5 & Cell 6& Cell 7& Cell 8&Cell 9& Cell 10 \\
{\bf conduct.}                  & \multicolumn{3}{c|}{} & \multicolumn{2}{c|}{} & \multicolumn{3}{c||}{} &   &   \\
$(\textrm{mS}/\textrm{cm}^{2})$ & \multicolumn{3}{c|}{} & \multicolumn{2}{c|}{} & \multicolumn{3}{c||}{} &   &   \\\hline
$g_{leak}$      	            & 0,1	& 0,0205& 0,0133& 0,15  & 0,038	 & 0,01  & 0,01  & 0,1	 &0,05  & 0,1\\
$g_{Na}$            	        & 50	& 56	& 10    & 	50  & 	58   & 	50   & 	50   & 	50   & 	3   & 35\\
$g_{K}$      	                & 5     & 6     & 21    & 10    & 3,9	 & 5     & 5     & 4,2   &  5   & 9\\
$g_{M}$                         & 0,07	& 0,075 & 0,098 & -  	& 0,0787 & 0,03	 & 0,03  & 0,042 & -    & -\\
$g_{L}$     	                & -     & -     & -     & -     & -      & 0,1   & 0,2   & 0,12  & -    & -\\
$g_{T}$                         & -     & -     & -  	& -     & -  	 & -     & -     & -  	 & 5    & 5\\\hline
{\bf Reversal}                  &  	    &       & 	    & 	    & 	 	 &  	 &  	 &       &      & \\
{\bf potential}                 &  	    &       & 	    & 	    & 	 	 &  	 &  	 &       &      & \\
(mV)                            &  	    &       & 	    & 	    & 	 	 &  	 &  	&        &    & \\\hline
$E_{leak}$	                    & -70	& -70,3	& -56,2	& -70	& -70,4	 & -70   & -70  & -75	 & -70& -65\\
$E_{Na}$ 	                    & 50	& 50	& 50	& 50	& 50     & 50    & 50	& 50	 & 50 & 55\\
$E_{K}$ 	                    & -90	& -90	& -90	& -90	& -90	 & -90   & -90	& -90	 & -90& -90\\
$E_{Ca}$        	            & -    	&  -	& - 	& -	    & -      & 120   & 120	& 120	 & -  & -\\
$E_{T}$          	            & -   	& -  	& - 	& -	    & -      & -     & -	& - 	 & 0  & -\\\hline
{\bf Other}                     &  	    &       & 	    & 	    & 	 	 &  	 &  	&        &    & \\
{\bf param.    }                &  	    &       & 	    & 	    & 	 	 &  	 &  	&        &    & \\\hline
$V_T$ (mV)          	        & -61,5	& -56,2 & -65,4 & -61,5 & -57.9  & -56,2 & -56,2& -58    & 	- & -\\
$V_x$ (mV)                      &  	-   &   -   & 	-   & 	 -  & 	 -	 &  -	 &  -	&   -    &  - & 5 \\
$\tau_{max}$ (mS)          	    & 4000	& 608   & 	934 & 	-   & 	502  & 4000  & 4000	& 1000   & 	- & -\\
$\phi$                          &  	-   &   -   & 	-   & 	 -  & 	 -	 &  -	 &  -	&   -    &  - & 5 \\
C ($\mu F$)          	        & 0,29	& 1     & 	1   & 	0,14& 	1    & 0,29  & 0,29	& 0,29   & 	1   & 1\\
\hline\hline

\end{tabular}
\end{center}

\begin{tabular}{ll}
&\\
Cell 1:& RS cell as observed from ferret Visual Cortex in vitro\\
Cell 2:& RS excitatory cell as observed from somatosensory cortex in vitro\\
Cell 3:& RS inhibitory cell as observed from somatosensory cortex in vitro\\
Cell 4:& FS cell as observed from ferret Visual Cortex in vitro\\
Cell 5:& FS cell as observed from somatosensory cortex in vitro\\
Cell 6:& IB cell as observed from guinea pig somatosensory cortex in vitro\\
&(Initial burst followed by adaptive action potentials)\\
Cell 7:& IB cell as observed from guinea pig somatosensory cortex in vitro\\
&(Repetitive bursting)\\
Cell 8:& IB cell as observed from cat visual cortex\\
Cell 9:& TCR cell as observed from Mouse thalamocortical relay neuron\\
Cell 10:& RHI cell as observed from Rat hippocampal interneuron\\
\end{tabular}

\caption{The maximal conductances and reversal potential values corresponding to each of the neuron models. Other parameters are also reported.}
\label{table2}
}
\end{table}
\subsubsection{Neocortical neurons}
In the models used to account for the different firing patterns characterizing cells in the neocortex, $I_{Na}=g_{N{a}}m^{3}h(V-E_{N{a}}) \textrm{ and } I_K=g_{K}n^4(V-E_K)$ are the sodium and potassium currents responsible for action potentials. $I_M=g_Mp(V-E_K)$ is a slow voltage-dependent potassium current responsible for spike-frequency adaptation, and $I_L=g_{L}q^2r(V-E_{Ca})$ is a calcium current to generate bursting.

The gating variables $m$, $h$, $n$, $q$ and $r$ obey the standard kinetic equation, $\dot{x}=\alpha_x(1-x)-\beta_x x$, where $x=m,h,n,q,r$ and $\alpha_x$ and $\beta_x$ are voltage-dependent variables. $m$ and $h$ are sodium channels activation and deactivation variables, $n$ is potassium channels activation variable, and $q$ and $r$ are calcium activation and inactivation variables.
The gating variable $p$ associated with the slow voltage-dependent potassium current is governed by the following equation: $\dot{p}=(p_\infty(V)-p)/\tau_p(V)$.

The RS model consists of the sodium and potassium currents responsible for action potentials ($I_{Na}$ and $I_K$) and a slow voltage-dependent potassium current responsible for spike-frequency adaptation ($I_M$). The IB model includes the same ionic currents as in the RS ones and an additional L-type calcium current ($I_L$). The simple FS model accounts only for the currents responsible of spike generation (i.e. $I_{Na}$ and $I_K$), and reproduces the spiking properties of FS cells as observed from ferret Visual Cortex in vitro. To capture the firing characteristics of FS cells as recorded from somatosensory cortex in vitro a slow potassium current ($I_M$) has been added to the simple FS model.

\subsubsection{Thalamocortical relay neuron}
The thalamocortical relay neuron model consists of currents responsible for generating spikes, $I_{Na}=g_{Na}m_{\infty}^3h(V-E_{Na})$ and $I_K=g_K(0.75(1-h))^4(V-E_K)$, as well as, low-threshold calcium current, $I_T=g_Tp_{\infty}^2r(V-E_T)$.
The dynamics of the gating variables $h$ and $r$ obey respectively the equations $\dot{h}=(h_{\infty}-h)/\tau_{h} \textrm{ and } \dot{r}=(r_{\infty})-r)/\tau_r$.
This model achieves a single spike activity consisting of trains of action potentials whose frequency depends on the strength of depolarization. When exposed to depolarizing current pulse of constant amplitude it results in a train of action potentials with no frequency adaptation.
%
\subsubsection{Hippocampal interneuron}
Finally, the hippocampal interneuron model obeys the same current-balance equation (Eq. (1)), and consists of leak current and the spike-generating $Na^+$ and $K^+$ voltage-dependent ion current ($I_{Na}=g_{N{a}}m_{\infty}^{3}h(V-E_{N{a}})$ and $I_K=g_Kn^4(V-E_K)$). For the transient sodium current, the activation variable $m$ is assumed fast and substituted by its steady-state function $m_{\infty}=\alpha_m/(\alpha_m+\beta_m)$. This model has the ability to reproduce repetitive spikes at high frequencies in response to a constant injected current. It has a small current threshold (the rheobase $I_{stim}\simeq 0.25 \mu A/cm^2$), and the firing rate is as high as 400Hz for $I_{stim}=20\mu A/cm^2.$

In this work, we use the conductance-based model by Pospischil et al. \cite{Pospischil2008} to reproduce the main features of the typical firing pattern of RS, FS and IB neurons as observed in different cells in the neocortex. The computational models used to account for the firing characteristics of thalamocortical relay neuron and the hippocampal fast-spiking interneurons follow respectively the works by \cite{Sohal2002,Guo2008,Wang1996}. The forms of the activation and inactivation functions describing the ion currents of all these models are reported in Table (\ref{table1}). The the values of biophysical parameters characterizing each of these models are shown in Table (\ref{table2}), and give rise to ten different cells. Cell 1 is a RS cell from ferret visual cortex, Cells 2 and 3 are respectively excitatory and inhibitory RS cells from a rat somatosensory cortex, Cells 4 and 5 are respectively FS cells from ferret visual cortex and rat somatosensory cortex, Cells 6 and 7 are both IB cells from guinea-pig somatosensory cortex (Cell 6 is characterized by an initial burst followed by adaptive action potentials, while Cell 7 gives rise to repetitive bursting.), Cell 8 is an IB cell from a cat visual cortex, Cell 9 is a mouse thalamocortical relay cell, and finally Cell 10 account for a rat hippocampal interneuron.

\begin{figure}
\begin{center}
\begin{tabular}{c}
\includegraphics[width=1\textwidth]{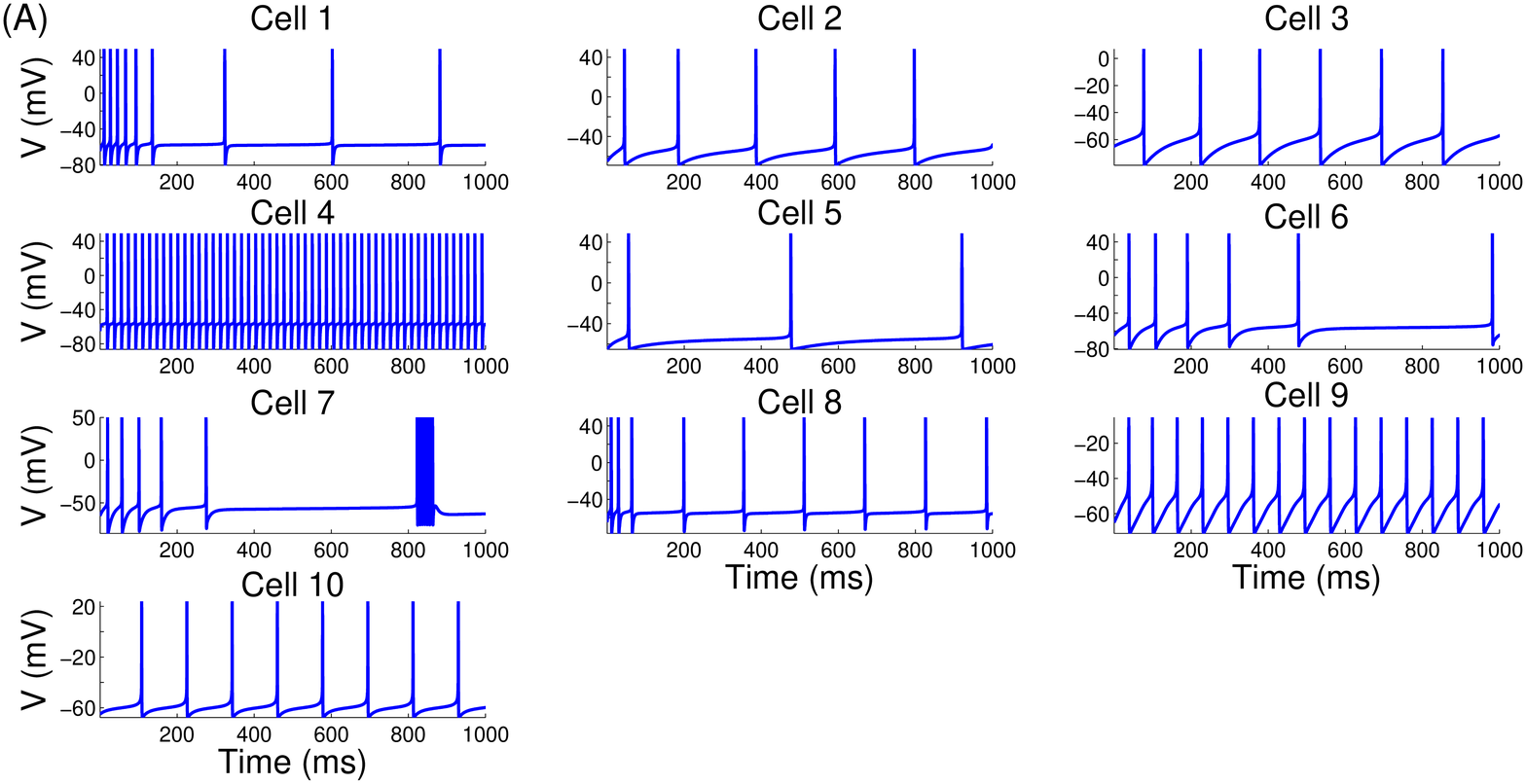}\\
\includegraphics[width=1\textwidth]{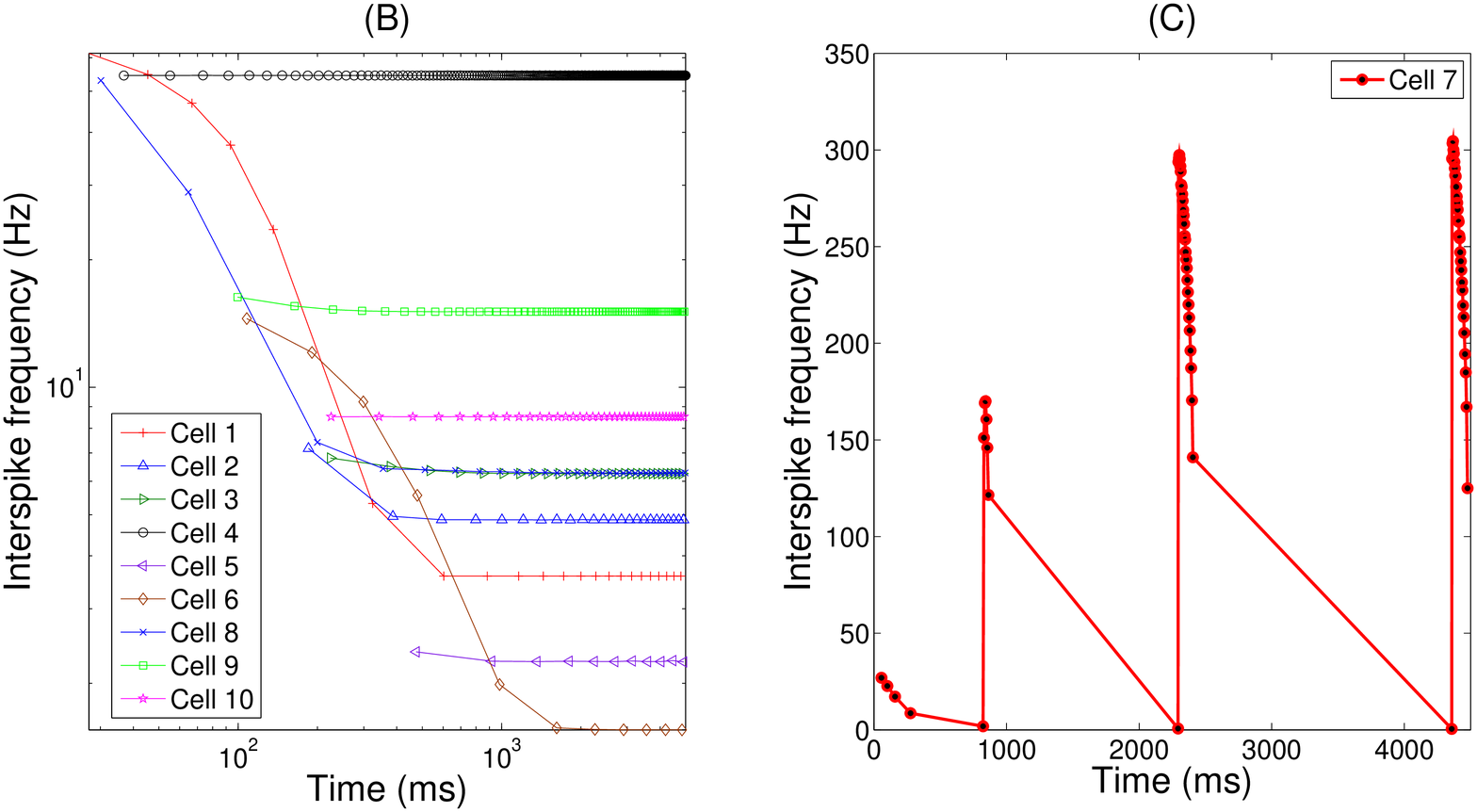}
\end{tabular}
\end{center}
\caption{ (A) Time course of membrane voltage from different spiking cells when stimulated by prolonged depolarizing stimulus slightly greater than threshold. The values of stimuli in $(\mu A/cm^2)$ are reported in Table 3.
(B) and (C) The interspike firing frequencies calculated from the interspike intervals as a function of time.
}
\label{fig1}
\end{figure}

\subsection{Ionic channels energy in spiking neuron models}

The procedure followed to find the electrochemical energy involved in the dynamics of the Hodgkin-Huxley circuit consisting of capacitor $C$ and three Na, K and L ionic channels has been reported in detail in \cite{Moujahid2011}. Here we extend this procedure to deduce the energy functions characterizing the dynamics of Hodgkin-Huxley-like models where additional ionic channels are present.
For the system given by Eq. (\ref{equ1}), the total electrical energy accumulated in the circuit at a given moment in time is,

\begin{equation}
H(t) = \frac{1}{2}\, CV^2+H_{l}+H_{N_{a}}+H_{K}+H_{M}+H_{L}+H_{T},
\label{equ2}
\end{equation}
where the first term in the summation gives the electrical energy accumulated in the capacitor and represents the energy needed to create the membrane potential $V$ of the neuron. The other six terms are the respective energies in the batteries needed to create the concentration jumps in chloride, sodium, potassium and calcium. The electrochemical energy
accumulated in the batteries is unknown. Nevertheless, the rate of electrical energy provided to the circuit by a battery is known to be the electrical current through the battery times its electromotive force.
Thus, the total derivative with respect to time of the above energy will be,
\begin{equation}
\dot{H}(t) = CV\dot{V}+I_{l}E_{l}+I_{N_{a}}E_{N_{a}}+I_{K}E_{K}+I_ME_{K}+I_LE_{Ca}+I_TE_{T}.
\label{equ3}
\end{equation}
where $E_{l}$, $E_{N{a}}$, $E_{K}$  and $E_{Ca}$ are the Nernst potentials of leakage, the sodium, potassium  and calcium ions in the resting state of the neuron. And $I_{l}$, $I_{N_{a}}$, $I_{K}$ ($I_{M}$), $I_{L}$ and $I_T$ are the ion currents of leakage, sodium, potassium, and calcium given by,
\begin{equation}
\begin{array}{ll}
I_{l}=g_{l}(V-E_{l}),\\
I_{N_a}=g_{N{a}}m^{3}h(V-E_{N{a}}),\\
I_{K}=g_{K}n^{4}(V-E_{K}), (\scriptstyle I_K=g_K(0.75(1-h))^4(V-E_K)  \textrm{ \scriptsize for TCR cell,})\\
I_M=g_Mp(V-E_K),\\
I_L=g_{L}q^2r(V-E_{Ca}),\\
I_T=g_Tp_{\infty}^2r(V-E_T)
\end{array}
\label{equ4}
\end{equation}

Substituting Eq. (\ref{equ1}) in Eq. (\ref{equ3}), we have for the energy rate in the circuit,

\begin{equation}
\begin{array}{ll}
\dot{H}(t) &= VI_{Stim}-I_l(V-E_l)-I_{N_{a}}(V-E_{N_{a}})-I_{K}(V-E_{K})\\
&-I_M(V-E_{K})-I_L(V-E_{Ca})-I_T(V-E_{T}).
\end{array}
\label{equ5}
\end{equation}

Finally, replacing in Eq. (\ref{equ5}) the ion currents by their expressions given by Eq. (\ref{equ4}), the electrochemical energy rate in each neuron will be given by,
\begin{equation}
\begin{array}{ll}
\dot{H}_{cell} = VI_{Stim}-\textrm{E}_{cell},\\
\end{array}
\label{equ6}
\end{equation}
where $\textrm{cell=RS, FS, IB, TCR, RHI}$ refers to the five conductance-based models considered in this work. This energy derivative provides the total derivative of the electrochemical energy in the neuron as a function of its state variables. The first term in the right hand summation represents the electrical power given to the neuron via the different junctions reaching the neuron and the $\textrm{E}_{cell}$ represents the energy per second consumed by the ion channels (ion channels energy) involved in the dynamics of each cell model. The forms of the ion channels energies are given as:
\begin{equation}
\begin{array}{ll}
\textrm{E}_{RS} = \scriptstyle g_{l}(V-E_{l})^2+g_{N_{a}}m^3h(V-E_{N_{a}})^2+g_{K}n^4(V-E_{K})^2+g_Mp(V-E_{K})^2,\\
\textrm{E}_{FS} = \scriptstyle g_{l}(V-E_{l})^2+g_{N_{a}}m^3h(V-E_{N_{a}})^2+g_{K}n^4(V-E_{K})^2,\\
\textrm{E}_{IB} = \scriptstyle g_{l}(V-E_{l})^2+g_{N_{a}}m^3h(V-E_{N_{a}})^2+g_{K}n^4(V-E_{K})^2+g_Mp(V-E_{K})^2+g_Lq^2r(V-E_{Ca})^2,\\
\textrm{E}_{TCR} = \scriptstyle g_{l}(V-E_{l})^2+g_{N_{a}}m^3h(V-E_{N_{a}})^2+g_{K}(0.75(1-h))(V-E_{K})^2+g_Tp_{\infty}^2r(V-E_T)^2,\\
\textrm{E}_{HFI} = \scriptstyle g_{l}(V-E_{l})^2+g_{N_{a}}m^3h(V-E_{N_{a}})^2+g_{K}n^4(V-E_{K})^2,
\end{array}
\label{equ7}
\end{equation}
These equations permit evaluation of the total energy consumed by a given neuron and also give information about the consumption associated to each of the ionic channels. Fundamentally, for each cell, this rate of energy expressed in nJ/s must be replenished by the ion pumps and metabolically supplied by hydrolysis of ATP molecules in order to maintain the cell's spiking activity.

\subsection{Ion-counting based energy consumption}
The energy consumption can be also computed based on the amount of $Na^+$ load ($Q_{N_a}$) crossing the membrane during a single action potential, which it can be estimated integrating the area under the total $Na^+$ current curve (see Fig. (2)) for the duration of stimulus presentation.

Since the $Na^+$ gradient is maintained primarily by the activity of the $Na^+/K^+$ ATPase which mediates influx of 2 K$^+$  in exchange for 3 Na$^+$  per one ATP molecule consumed, the number of ATP moles ($ATPmols$) can be computed as  $\frac{Q_{Na}}{3\textrm{eNA}}$, where e=1.602 $10^{-19}$ C is the electric charge and NA=6.022 $10^{23}$ is the Avogadro constant \cite{Attwell2001}. For the other hand, the energy available for chemical work from the ATP concentration is measured by the free energy of ATP hydrolysis ($F_{ATP}$). This allows an estimate of the metabolic energy associated with ionic pumping as the free energy times the number of ATP moles, that is, $F_{ATP}*ATPmols$.

The value of the free energy depends to some extent on the internal chemical state of the cell, and it has been reported to be in the range from 46 kJ/mol to 62 kJ/mol of ATP \cite{Jansen2003}, but in this work, we set $F_{ATP}=50 \textrm{ kJ/mol}$.

Moreover,  the integral of the area under the instantaneous energy function for a given cell (Eqs. (\ref{equ7})) provides an estimate of the neuron energy, which is used to computed the efficiency of the ATP hydrolysis as the ratio of that energy consumption to the number of ATP moles, that is, $\frac{1}{ATPmols}\int E_{cell}$. This ratio expressed is kJ/mol gives us the opportunity to check if our method of calculation of the actual energy consumption by the pump
and the number of ATP molecules involved are consistent with other data in the literature. The calculated values of this ratio are reported in figures (3) and (6) as hydrolysis of ATP molecules.

\begin{figure}
\begin{center}
\includegraphics[width=1\textwidth]{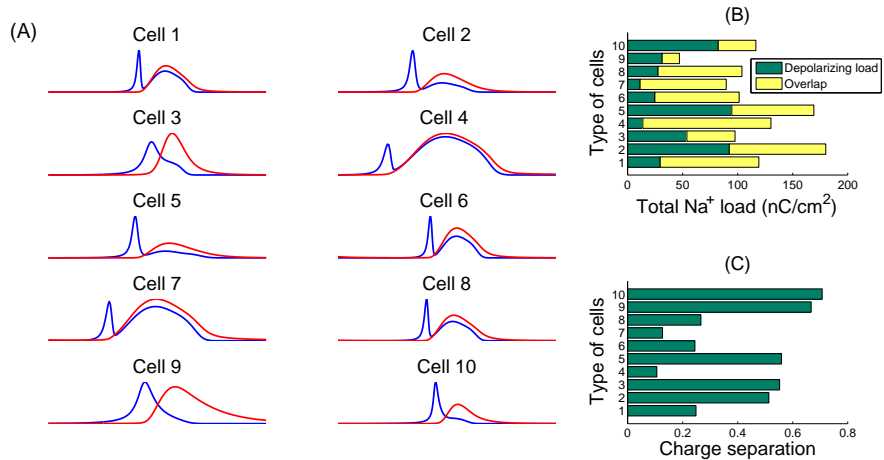}
\end{center}
\caption{ (A) The sodium (blue line) and potassium (red line) currents during an action potential showing different degree of overlap. Sodium current is reversed for comparison. The neutralized flux is measured as the difference between the total Na$^+$ load and the Na$^+$ minimum charge transfer necessary for the depolarization of the action potential.
(B) The total Na$^+$ load per action potential for each cell depicted as the sum of the capacitive minimum and the overlap loads. (C) The charge separation computed as the ratio of the capacitive minimum to the total Na$^+$ load. This ratio show how efficiently are the considered cells in generating action potentials. Action potentials have been generated for the values of stimuli reported in Table 3.
}
\label{fig2}
\end{figure}

\begin{figure}
\begin{center}
\includegraphics[width=1\textwidth]{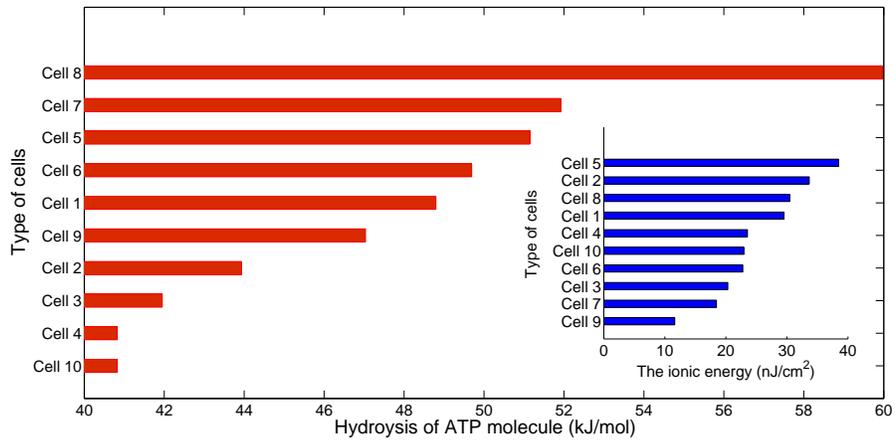}
\end{center}
\caption{ The free energy of ATP hydrolysis in kJ/mol computed as the ratio of the ionic energy and the number in mole of ATP molecules per membrane unit area. The values of the free energy liberated range from 40.82 kJ/mol in the case of the rat hippocampal interneuron to 59.95 kJ/mol in the intrinsically bursting cell from a cat visual cortex.  The inset displays the corresponding values of ionic energy in nJ per cm$^2$ calculated from Eq. (7).
}
\label{fig3}
\end{figure}

\begin{table}
{\scriptsize
\begin{tabular}{l||ccc|cc|ccc||c||c}
\hline\hline
                                & \multicolumn{3}{c|}{\bf RS cells} & \multicolumn{2}{c|}{\bf FS cells} & \multicolumn{3}{c||}{\bf IB cells} & {\bf TCR}  & {\bf RHI}  \\
                                &  	    &       & 	    & 	    & 	 	 &  	 &  	 &       &      & \\ \hline
             	                & Cell 1& Cell 2& Cell 3& Cell 4& Cell 5 & Cell 6& Cell 7& Cell 8&Cell 9& Cell 10 \\
Frequency (Hz)	                    & 5 	& 5 	& 6 	& 54    & 2 	 & 2     & 15 	 & 7	 & 15   & 9\\\hline
$Na^+$  load 	                & 174   & 207   & 134   & 162   & 217	 & 132   & 103   & 147   & 69   & 163\\
$K^+$ load    	                & 141   & 214   & 150   & 156   & 197    & 137   & 117   & 133    & 79   & 127\\
Capacitive minimum              & 65	& 108   & 70    & 22 	& 129    & 37	 & 15    & 51    & 55   & 125\\
Overlap load                    & 109   & 99    & 64	& 140   & 88	 & 95    & 88    & 96 	  & 14   & 38\\\
$(\textrm{nC}/\textrm{cm}^{2})$ &  	    &       & 	    & 	    & 	 	 &  	 &  	 &        &      & \\\hline
Charge separation               &0.38   &0.52   &0.52 	&0.14   &0.60    &0.28   &0.14   &0.35   &0.79 	&0.77\\\hline
ATP Pmole            	        & 0.60	& 0.72	& 0.46  & 	0.56& 	0.75 & 	0.46 & 	0.36 & 	0.51 & 	0.24& 	0.56\\\hline
Metabolic Energy	            & 30	& 36	& 23	& 28	& 38	 & 23    & 18    & 25    & 12   & 28  \\
(ion-counting method)           &  	    &       & 	    & 	    & 	 	 &  	 &  	&        &      & \\
$(\textrm{nJ}/\textrm{cm}^{2})$ &  	    &       & 	    & 	    & 	 	 &  	 &  	 &        &      &  \\\hline
Ionic Energy	                & 30    & 34	& 20 	& 24	& 38 	 & 23    & 18   & 30	 & 12   & 23 \\
(From Eq. (\ref{equ7}))         &  	    &       & 	    & 	    & 	 	 &  	 &  	&        &      & \\
$(\textrm{nJ}/\textrm{cm}^{2})$ &  	    &       & 	    & 	    & 	 	 &  	 &  	 &        &      &  \\\hline
ATP Hydrolysis                  &49.14	&47.03	&43.93	&41.96	&51.15	 &49.70	&51.91	&59.95	&48.78	&40.82\\
 (kJ/mol)	 	                &  	    &       & 	    & 	    & 	 	 &  	 &  	&        &      & \\\hline
Stimulus                        & 1.4	& 0.7   & 	0.15& 	1.75& 	0.8  & 0.25  & 0.25	& 2.25   & 	0.44 & 	0.20\\
$(\mu A/cm^2)$	                &  	    &       & 	    & 	    & 	 	 &  	 &  	&        &      & \\\hline\hline

\end{tabular}

\caption{Ionic flux and energy demands of single action potentials from different spiking cells when stimulated by prolonged current stimulus slightly greater than threshold. The metabolic energy refers to the energy computed according to the ion-counting approach (Section 2.3), and the ionic energy accounts for the electrochemical energy computed as the integral of the energy functions given by Eq. (\ref{equ7})  }
\label{table3}
}
\end{table}

\section{Results}

\subsection{Firing characteristics of action potentials}
The different possible combinations of conductances in the neuron models considered in this work give origin to a variety of action potentials with different waveforms and firing characteristics. Figure (\ref{fig1}) reports, in part (A), trains of action potentials from different neuron models when exposed to prolonged stimuli of constant magnitude slightly greater than threshold (See Table 3 for a description of the current stimuli values). The inter-spike firing frequencies calculated from the inter-spike time intervals are reported in parts (B) and (C). As it can be seen, of all regular-spiking cells, the RS cell from ferret visual cortex (Cell 1) exhibits the most pronounced adaptation with an inter-spike firing frequency declining from about 62 to 3.5 Hz.  Likewise, the IB cell as recorded respectively from guinea pig somatosensory cortex (Cell 6) and cat visual cortex (Cell 8), initially generate action potentials at high firing frequency (of about 55 Hz) that decreases rapidly to low values within a short time span. Meanwhile, the repetitive IB cell from guinea pig somatosensory cortex (Cell 7) shows repetitive intrinsic bursting with a mean inter-burst frequency of about 1 Hz, and an inter-spike frequency ranging from 300 Hz to 150 Hz. The others cell models (i.e. Cells 4, 5, 9 and 10) respond to depolarizing stimuli by generating action potentials with little or no adaptation.

\subsection{Efficient action potentials and overlapping}

Figure (\ref{fig2}) shows in part (A) the shape of the sodium (in blue line) and potassium currents corresponding to a particular action potential of the previously described trains. The sodium current (in blue line) is negative but has been depicted with a positive sign for a better appreciation of the great extent of its overlapping with the potassium current (red line). Note that as sodium and potassium currents are both of positive charges but moving in opposite directions of the cell's membrane they neutralize each other to the extent of their mutual overlap. The sodium charge that is not counterbalanced by simultaneously flowing potassium charge is much smaller for a greater overlap. This unbalanced load corresponds to the minimum charge ($Q_{min}$) needed for the depolarization of the action potential. The sum of this capacitive minimum and the overlap load gives the total Na$^+$ loads per action potential. Recollected values are depicted in Fig. (\ref{fig2})(B).

To quantify the sodium and potassium currents overlap we calculated the dimensionless charge separation as the ratio of the capacitive minimum and the total Na$^+$ charge per action potential. That is, $Q_{Separation}=Q_{min}/Q_{Na}$, this measure allows quantifying how efficiently sodium flux is used for action potential depolarization \cite{Alle2009}. In view of the above, it seems that rat hippocampal interneuron (Cell 10) and mouse thalamocortical relay neuron (Cell 9) are the most efficient in generating action potentials with a charge separation approaching 80\%, which shows that most of sodium load are confined to the rising phase of the action potential. For the other hand, the FS cell from a ferret Visual Cortex (Cell 4) and the repetitive IB cell from a guinea pig somatosensory cortex in vitro (Cell 7) are the most inefficient with only 14\% of sodium entry used for spike depolarization. In fact, the sodium influx during an action potential in thalamocortical relay neuron was 69 nC/cm$^2$, consuming 0.24 pmol/cm$^2$ of ATP molecules per spike, while the FS cell from somatosensory cortex moves 217 nC/cm$^2$ of sodium ions per spike consuming 0.75 pmol/cm$^2$ of ATP molecules. Other models including RS cells from somatosensory cortex (Cell 2) and ferret visual cortex (Cell 1) show high values of Na$^+$ loads of about 207 nC/cm$^2$ (0.72 ATP pmol/cm$^2$) and 174 nC/cm$^2$ (0.60 ATP pmol/cm$^2$) respectively. These estimates agree with other data reported in literature \cite{Sengupta2010,Carter2009,Alle2009}.

\subsection{Energy}
Assuming that one ATP molecule hydrolyzed under normal physiological conditions liberates of about 50kJ/mol, our calculations of the energy demands in nJ per cm$^2$ based on the ion-counting method give energy consumptions ranging from 12 nJ/cm$^2$ for a thalamocortical relay cell (Cell 9) to 38 nJ/cm$^2$ for a FS cell from somatosensory cortex (Cell 5).

For the other hand, following the ion-channels energy functions described in Eq. (7), the calculated values of the energy consumption give results which are in excellent agreement with those computed according to ion-counting approach. Indeed, according to the Wilcoxon rank-sum test \cite{Gibbons2011} performed on the data reported in Table 3, the difference between the values of metabolic and ionic energy is statistically not significant at the 5\% significance level with a p-value=0.8194.

The values of the energy consumption needed for the restoration of concentration gradients after an action potential was determined by integrating over long period of time the area under the instantaneous ion channels energy curve divided by the number of spikes, which gives the energy consumption of a single spike. In Figure (\ref{fig3}) we report the free energy of ATP hydrolysis computed as the ration of the ionic energy and the number of ATP molecules. This ratio is expressed in kJ/mol and gives values ranging from 40.82 kJ/mol in the case of a rat hippocampal interneuron (Cell 10) and value of about 60 kJ/mol for an IB cells as observed from cat visual cortex (Cell 8). These values are in nice agreement with the values of the free energy of ATP hydrolysis reported in the literature \cite{Jansen2003,Ereciriska1989}.

\begin{figure}
\begin{center}
\begin{tabular}{c}
\includegraphics[width=1\textwidth]{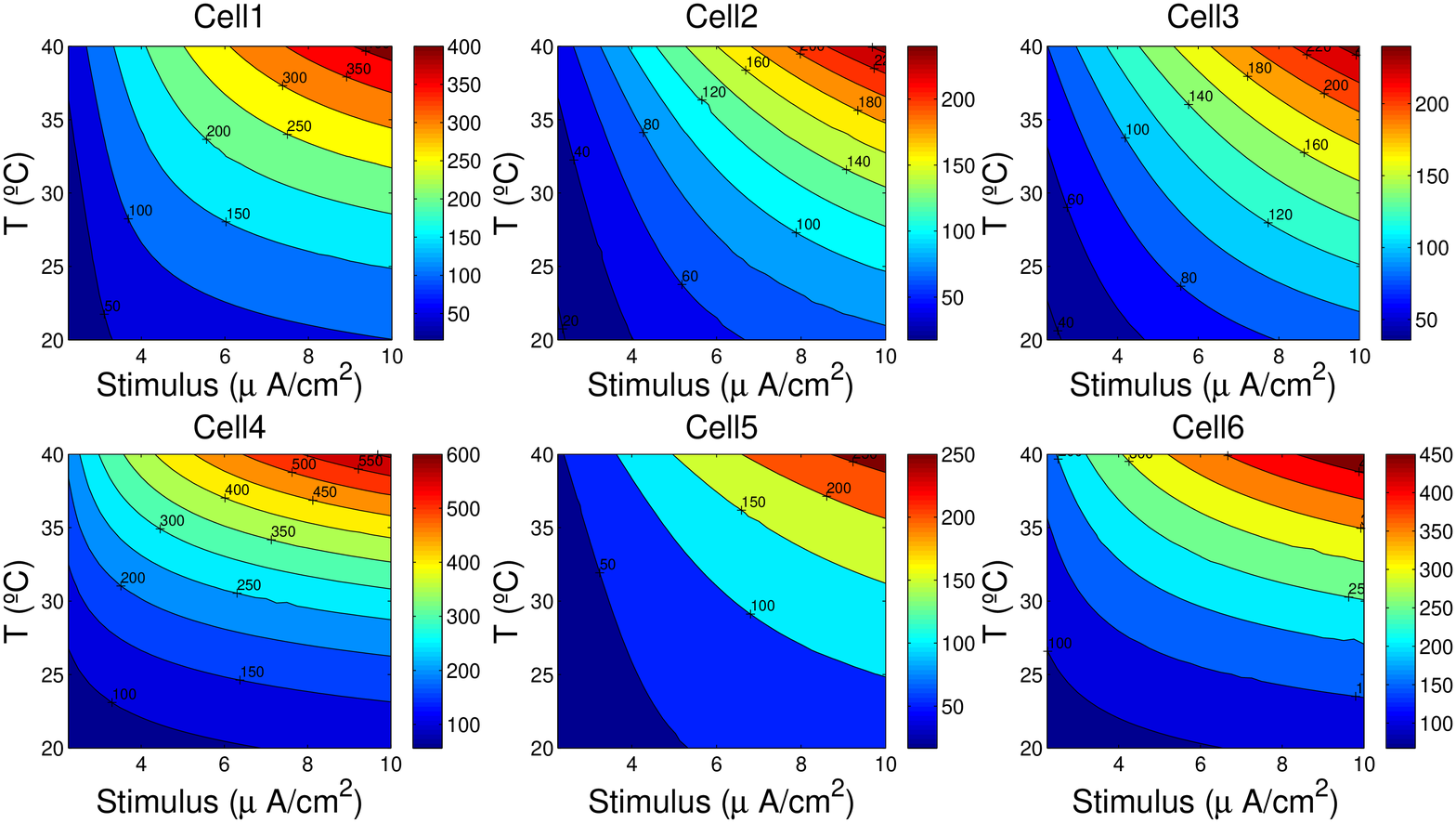}\\
\includegraphics[width=1\textwidth]{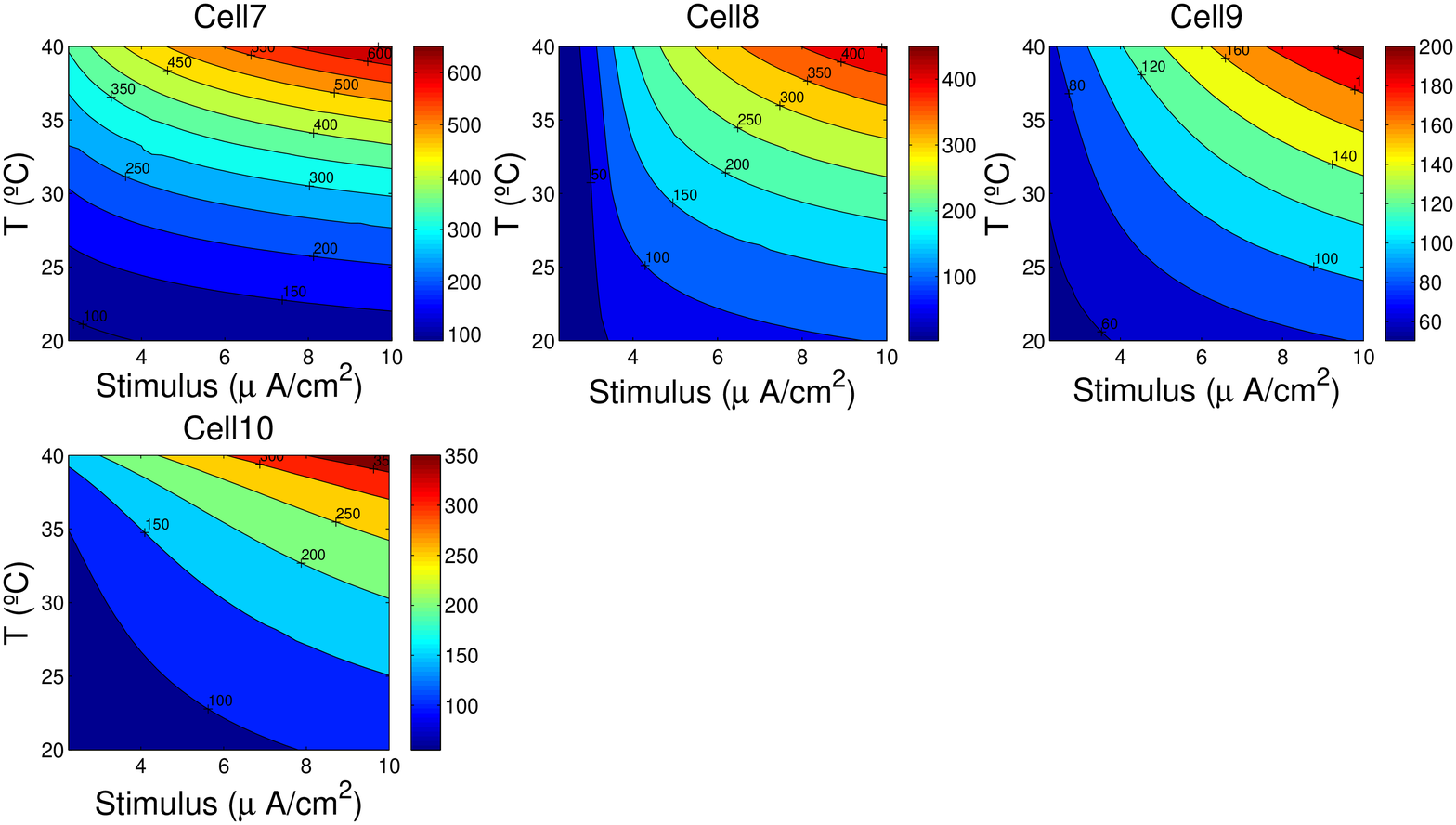}
\end{tabular}
\end{center}
\caption{The firing frequency (Hz) as a function of the cells temperature and the external depolarizing stimulus. We consider a range of temperature between 20 $^\circ$C and 40 $^\circ$C, and a stimulus varying from 2.25 to 10 $\mu A/cm^2$. }
\label{fig4}
\end{figure}

\subsection{Temperature and stimulus intensity}
Because neural properties are temperature dependent, we have analyzed in more detail the values of energy consumption and the free energy of ATP hydrolysis as a function of temperature and the magnitude of the external stimuli. Varying the cell temperature in the neuron models was down multiplying the rates of change of the activation and inactivation gating variables by a factor $k = 2.78^{(Temp-36)/10}$, and considering a range of temperature between 20 and 40 $^\circ$C. The stimuli were varied from 2.25 to 10 $\mu A/cm^2$.

\begin{figure}
\begin{center}
\begin{tabular}{c}

\includegraphics[width=1\textwidth]{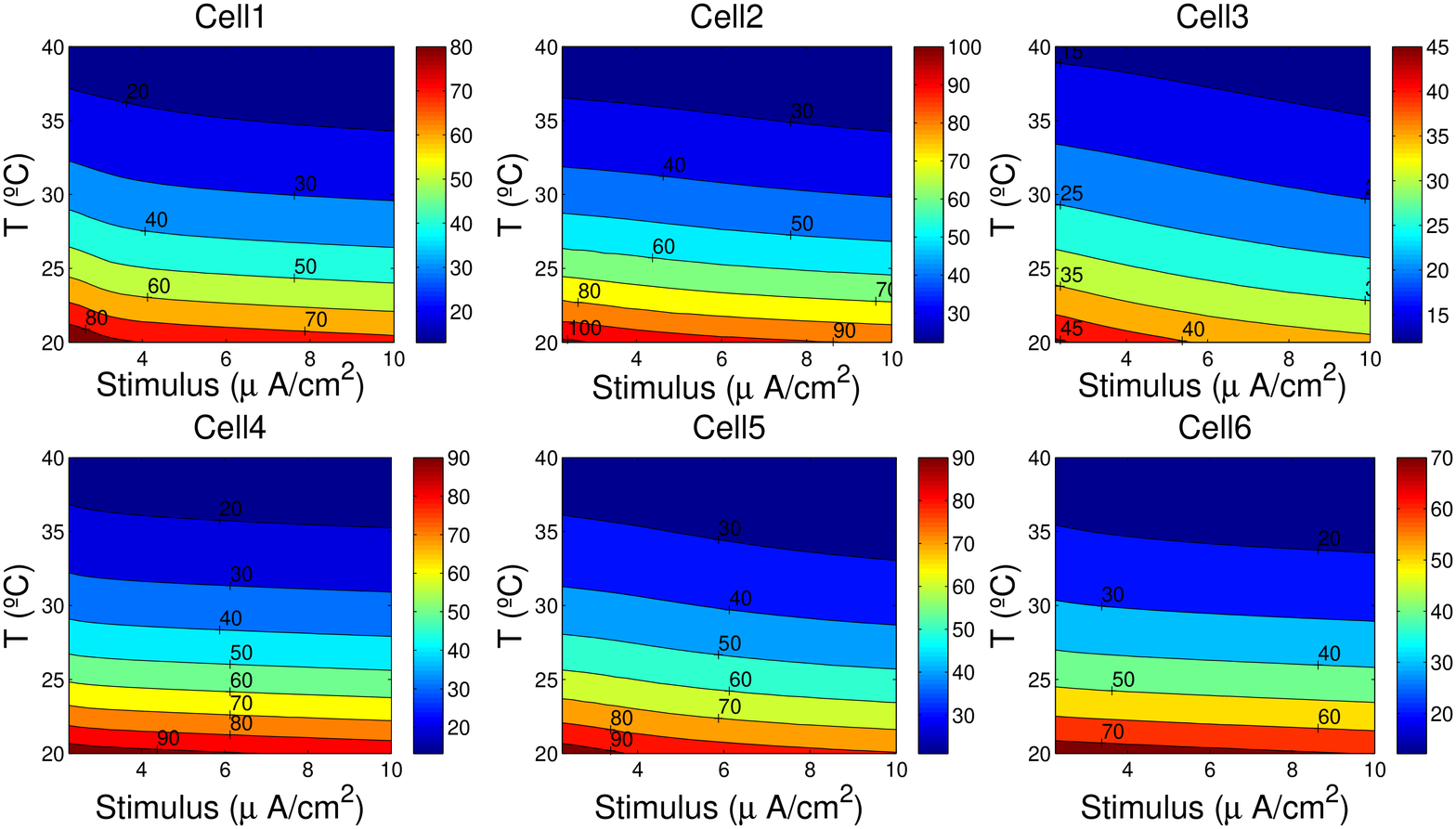}\\
\includegraphics[width=1\textwidth]{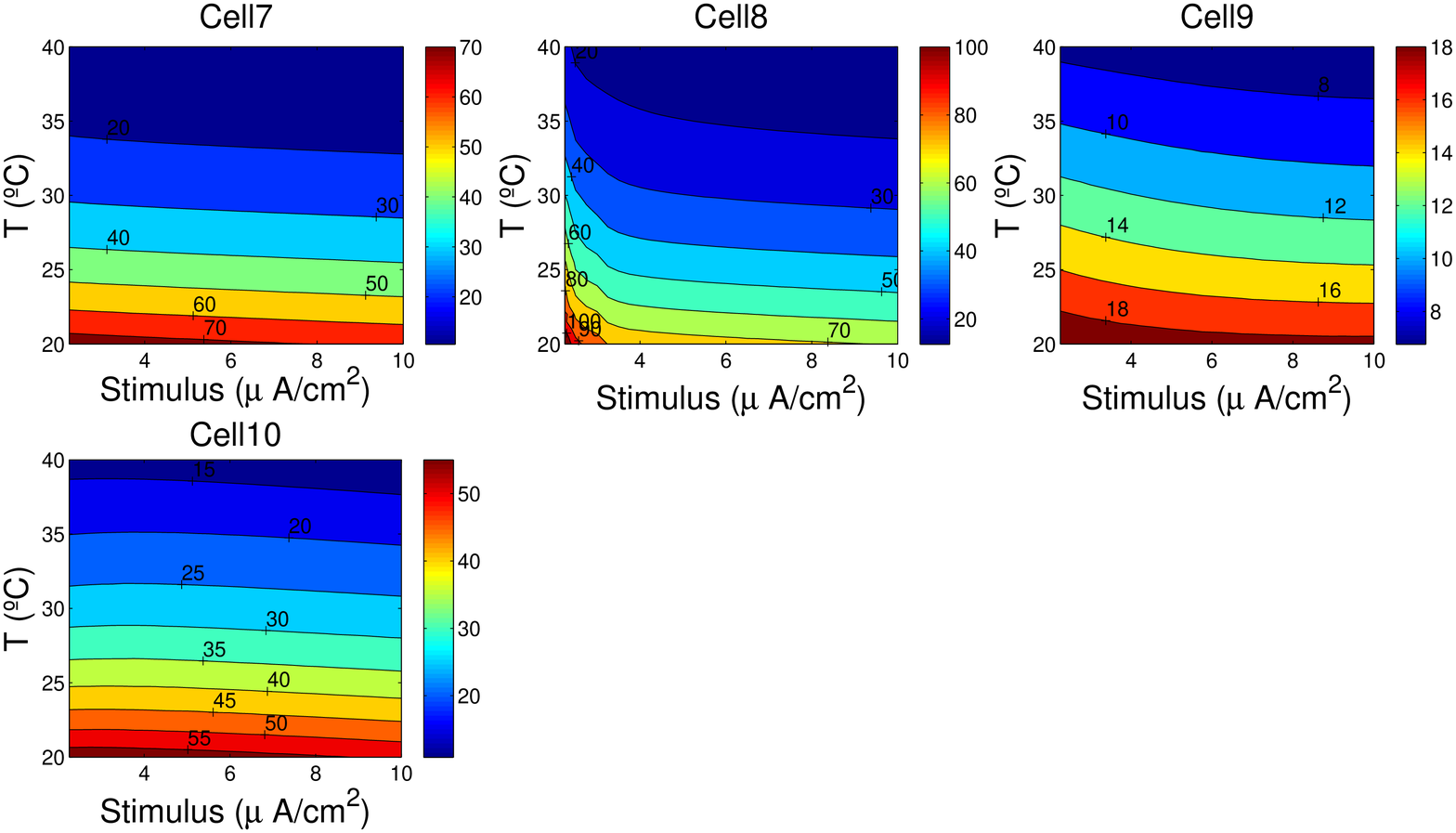}
\end{tabular}
\end{center}
\caption{The ionic energy $(\textrm{nJ}/\textrm{cm}^{2})$ as a function of the cells temperature and the external depolarizing stimulus. We consider a range of temperature between 20 $^\circ$C and 40 $^\circ$C, and a stimulus varying from 2.25 to 10 $\mu A/cm^2$. }
\label{fig5}
\end{figure}

The contour plots in Figure (\ref{fig5}) show, for each cell, how the energy consumption in nJ/cm$^2$ per spike depends upon cell temperature and input stimulus. As it can be seen the lowest energy consumption values are achieved for higher temperatures. In fact, the increased firing frequencies induced by higher temperatures and stimulus (see Figure (\ref{fig4})) imply more efficient use of sodium entry due mainly to the reduced overlap load between inward Na$^+$ current and outward K$^+$ current \cite{Moujahid2012}. Indeed, for all cells, the overlap load undergoes a reduction of about 90\% with an increase of cell temperature in a wide range of stimulus intensity. This efficient use of sodium entry is usually accompanied by a reduction of the energy demands needed for the restoration of concentration gradients.
For example, the rat hippocampal interneuron (Cell 10) at 20 $^\circ$C and an input stimulus of 2.25 $\mu$ A/cm$^2$ fires spikes at 55 Hz and consumes of about 58 nJ/cm$^2$.
While at a higher temperature of about 40 $^\circ$C for a wide range of stimulus intensity, the energy consumption is about 5-fold lower that the energy needed at low temperatures. This 5-fold decrease in energy consumption corresponds to 7-fold increase in firing frequency for an input stimulus of 10 $\mu$ A/cm$^2$. Likewise, the IB cell as observed from cat visual cortex (Cell 8) consumes of about 109 nJ/cm$^2$ at 20 $^\circ$C and lower stimulus, however at 40 $^\circ$C it requires only 13 nJ/cm$^2$ which represents of about 8-fold decrease in its energy demand. At this temperature the values of energy consumption are maintained for a wide range of input stimulus.

\subsection{The free energy of ATP hydrolysis}

\begin{figure}
\begin{center}
\begin{tabular}{c}
\includegraphics[width=1\textwidth]{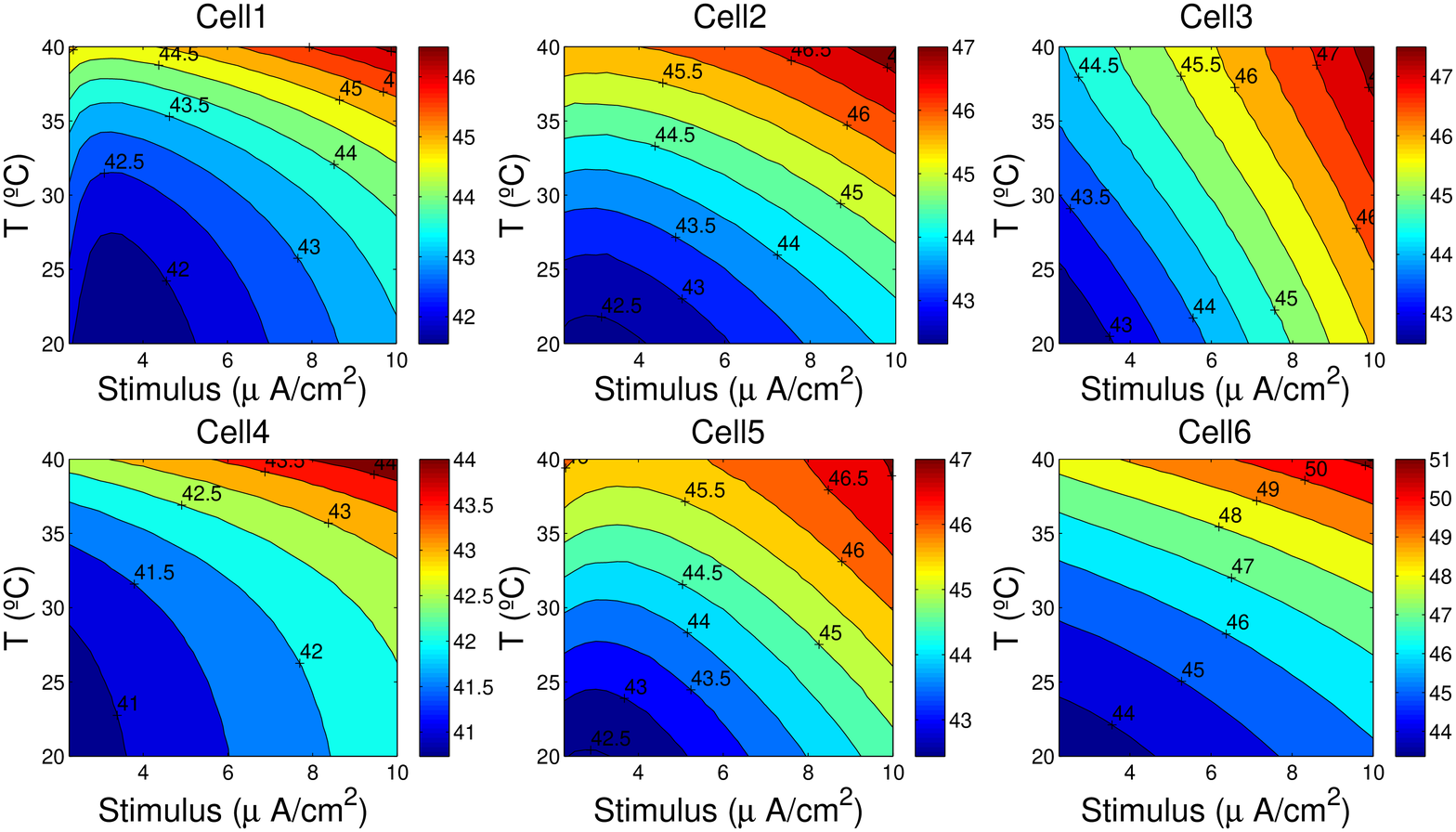}\\
\includegraphics[width=1\textwidth]{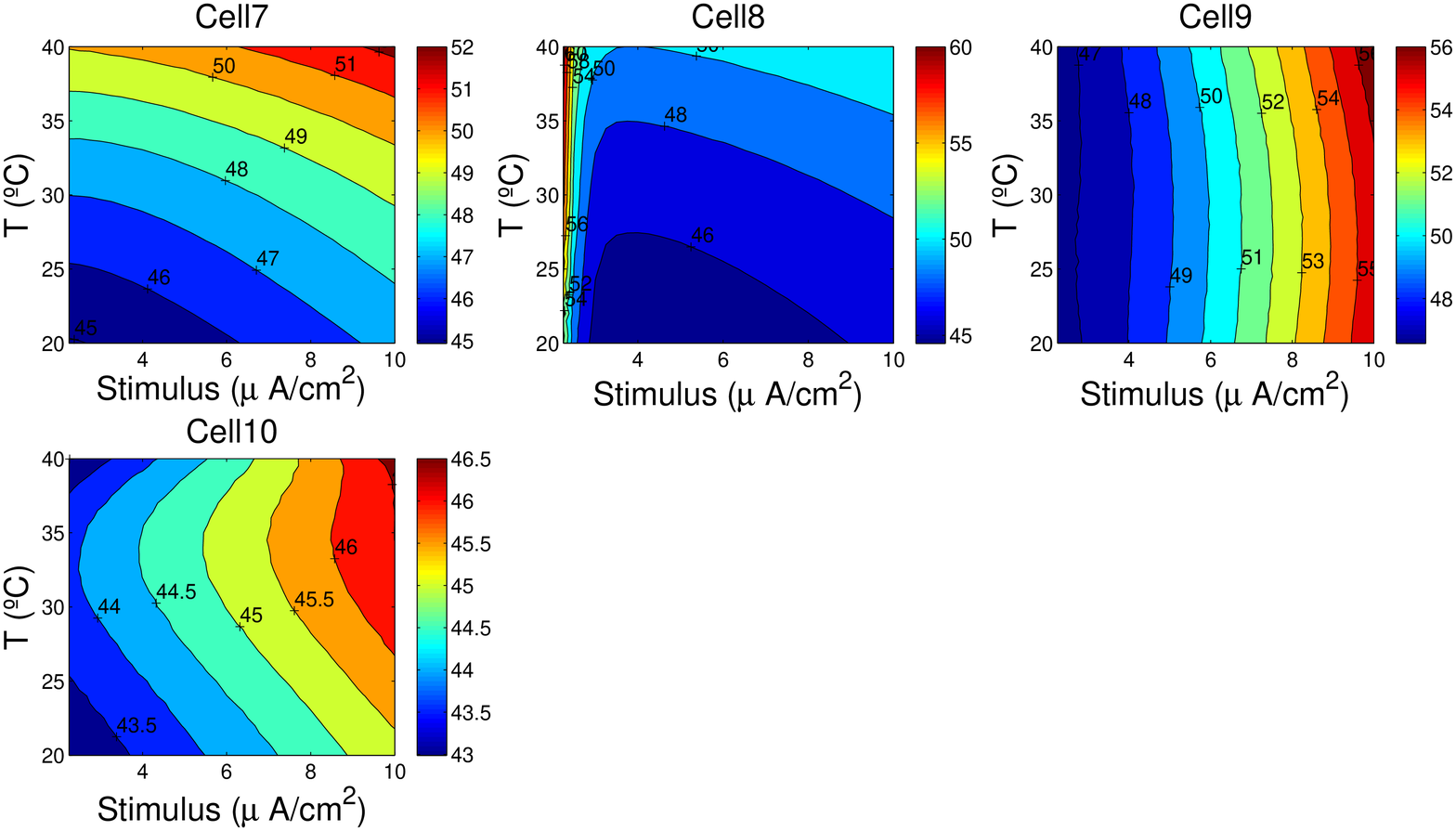}
\end{tabular}
\end{center}
\caption{Hydrolysis of ATP molecule in kJ/mol as a function of the cells temperature and the external depolarizing stimulus. We consider a range of temperature between 20 $^\circ$C and 40 $^\circ$C, and a stimulus varying from 2.25 to 10 $\mu A/cm^2$.}
\label{fig6}
\end{figure}

Figure (\ref{fig6}) shows the calculated values of ATP hydrolysis in kJ/mol as a function of temperature and stimulus intensity. In this case, we observe different patterns depending on the cells type. The RS and FS cells at higher temperatures and high input stimuli liberate respectively of about 12\% and 10 \% more free energy than at lower temperatures. For these cells the values of free energy liberated range from 40.73 kJ/mol in the case of FS cell from a ferret Visual Cortex (Cell 4) to 47.74 kJ/mol in an inhibitory RS cell from somatosensory cortex in vitro (Cell 3). The IB cells in guinea pig somatosensory cortex (Cells 6 and 7) show similar pattern as observed in regular and fast spiking cells, but in this case, increasing the cell temperature and stimulus intensity cause the free energy to experience an increase of more than 17\%  reaching values of about 52.27 kJ/mol. But, the more sensible cells to temperature and stimulus changes are IB cells as observed from cat visual cortex (Cell 8) and mouse thalamocortical relay neuron (Cell 9). Cell 8 shows values of free energy that range from 44.58 kJ/mol to 61.21 kJ/mol representing an increase of about 37\%. For this cell, the higher values of free energy are achieved for higher temperatures and low values of stimulus. In fact, the more pronounced pattern was observed for values of stimulus below 3 $\mu A/cm^2$, and a further increase of stimulus has seems to have little effect on free energy values. Conversely, Cell 9 accounting for a mouse thalamocortical relay neuron appears to be more sensitive to stimulus intensity than to temperature values. Increasing stimulus intensity from 2.25 to 10 $\mu A/cm^2$ moves the cell to liberate of about 22\% more free energy of ATP hydrolysis. Finally, for Cell 10 corresponding to a rat hippocampal interneuron, we can see that for values of input stimulus below 8 $\mu A/cm^2$, the free energy shows an hyperbole surface with a maximum around 33 $^\circ$C.

\section{Discussions}

Based on three computational neuron models (see Table 1), we used different combinations of voltage-dependent conductances to reproduce a variety of action potentials with different waveforms and firing characteristics as recorded from different cells in the neocortex, thalamus and hippocampus. Initially, each cell was stimulated by prolonged depolarizing stimulus slightly greater than threshold. We calculated in nC/cm$^2$ the amount of Na$^+$ load crossing the membrane during a single action potential and the corresponding metabolic energy nJ/cm$^2$ needed for the restoration of concentration gradients after an action potential following two different approaches (see Section 2). According to our calculations, for all the cells considered in this work, these two approaches have led to similar values of energy demands which confirms the consistency of our method to estimate the metabolic consumption. In fact, the range of values of our estimations of the free energy of ATP hydrolysis (see Figures 3 and 6) agree with other data reported in the literature \cite{Jansen2003,Ereciriska1989}.

An initial classification of the different cells (see Table 3, Figure 3) places the thalamocortical relay cell (Cell 9) as the most energy-efficient cell consuming 12 nJ/cm$^2$ for each spike generated and liberating a free energy of ATP hydrolysis of about 48.78 kJ/mol. For the other hand, the FS cell from somatosensory cortex (Cell 5) is the least energy-efficient requiring 38 nJ/cm$^2$ per spike and liberating 51.15 kJ/mol of free energy from ATP hydrolysis.

To quantify how the energy consumption in nJ/cm$^2$ and then the free energy of ATP hydrolysis depend upon cell temperature and input stimulus, we analyzed in more detail the values of these measures for a wide range of cell temperatures and stimulus intensities. The results give a clear evidence that the metabolic energy associated to the generation of action potentials is higher at lower temperature. Indeed, the spiking activity of the cells at higher temperatures involves less ion flux during the generation of an action potential. In all the neurons studied in this work, increasing the cell temperature from 20 to 40 $^\circ$C resulted in decreases by about 80\% in both the sodium and potassium flux. Obviously, this decrease in ion flux should be accompanied by a decrease in energy consumption. In fact, results reported in Figure 5 support this conjecture showing that low values of energy are achieved at higher temperatures. For a given temperature and stimulus intensity, for example, at 36 $^\circ$C and a stimulus value of 7 $\mu A/cm^2$, the energy consumptions range from 8.42 nJ/cm$^2$ in the case of cell 9 to 26.8 nJ/cm$^2$ in the case of cell 5 and 28.5 nJ/cm$^2$ in cell 2. The other cells show value of energy between 15 and 19 nJ/cm$^2$. While at 40 $^\circ$C and maintaining the same input stimulus (7 $\mu A/cm^2$) the values of energy consumptions reported above decreased by 17\%.

The approach reported in this work, contrary to ion-counting-based methods \cite{Attwell2001,Lennie2003}, does not require ion counting for estimating the metabolic energy consumption of the generation of action potentials, and gives us the opportunity to check in neuron models described by Hodgkin-Huxley type kinetics which ion counting gives the correct metabolic energy consumption.

On the other hand, since in many areas of the brain neurons are organized in populations of units with similar properties, it should be of interest to know about the metabolic cost of information processing of large population of interconnected neurons. For this end, we need to combine experimental studies of nervous systems with numerical simulation of large-scale brain models. Our approach to find an energy function that quantifies the physical energy associated to the states of a generic model neuron described by differential equations could be of great interest when studying population of coupled neurons.

%
%
%

%
%


\newpage

\begin{thebibliography}{39}
\expandafter\ifx\csname natexlab\endcsname\relax\def\natexlab#1{#1}\fi
\expandafter\ifx\csname bibnamefont\endcsname\relax
\def\bibnamefont#1{#1}\fi
\expandafter\ifx\csname bibfnamefont\endcsname\relax
\def\bibfnamefont#1{#1}\fi
\expandafter\ifx\csname citenamefont\endcsname\relax
\def\citenamefont#1{#1}\fi
\expandafter\ifx\csname url\endcsname\relax
\def\url#1{\texttt{#1}}\fi
\expandafter\ifx\csname urlprefix\endcsname\relax\def\urlprefix{URL }\fi
\providecommand{\bibinfo}[2]{#2}
\providecommand{\eprint}[2][]{\url{#2}}


\bibitem[{\citenamefont{Cauli et~al.}(1997)\citenamefont{Cauli, Audinat,
Lambolez, Angulo, Ropert, Tsuzuki, Hestrin and Rossier}}]{Cauli1997}
\bibinfo{author}{\bibnamefont{Cauli}, \bibfnamefont{B.}},
\bibinfo{author}{\bibnamefont{Audinat}, \bibfnamefont{E.}},
\bibinfo{author}{\bibnamefont{Lambolez}, \bibfnamefont{B.}},
\bibinfo{author}{\bibnamefont{Angulo}, \bibfnamefont{M. C.}},
\bibinfo{author}{\bibnamefont{Ropert}, \bibfnamefont{N.}},
\bibinfo{author}{\bibnamefont{Tsuzuki}, \bibfnamefont{K.}},
\bibinfo{author}{\bibnamefont{Hestrin}, \bibfnamefont{S.}}, \bibnamefont{and}
\bibinfo{author}{\bibnamefont{Rossier}, \bibfnamefont{J.}}, (\bibinfo{year}{1997}).
\emph{\bibinfo{title}{Molecular and physiological diversity of
cortical non pyramidal cells.}}
\bibinfo{journal}{J. Neurosci.} \textbf{\bibinfo{volume}{17}},
\bibinfo{pages}{3894–3906}.

\bibitem[{\citenamefont{Peters et~al.}(1984)\citenamefont{Peters and Jones}}]{Peters1984}
\bibinfo{author}{\bibnamefont{Peters}, \bibfnamefont{A.}}, \bibnamefont{and}
\bibinfo{author}{\bibnamefont{Jones}, \bibfnamefont{E. G.}},
(\bibinfo{year}{1984}).
\emph{\bibinfo{title}{Components of the Cerebral Cortex.}}
\bibinfo{Book}{Plenum, New York, 1984}.

\bibitem[{\citenamefont{Toledo-Rodriguez et~al.}(2003)\citenamefont{Toledo-Rodriguez,
Gupta, Wang, Wu and Markram}}]{Toledo-Rodriguez2003}
\bibinfo{author}{\bibnamefont{Toledo-Rodriguez}, \bibfnamefont{M.}},
\bibinfo{author}{\bibnamefont{Gupta}, \bibfnamefont{M.}},
\bibinfo{author}{\bibnamefont{Wang}, \bibfnamefont{Y.}},
\bibinfo{author}{\bibnamefont{Wu}, \bibfnamefont{C. Z.}}, \bibnamefont{and}
\bibinfo{author}{\bibnamefont{Markram}, \bibfnamefont{H.}}, (\bibinfo{year}{2003}).
\emph{\bibinfo{title}{in The Handbook of Brain Theory and Neural Networks. (ed. Arbib, M. A.) 719–725.}}
\bibinfo{Book}{MIT Press, Cambridge, Massachusetts., 2003}.

\bibitem[{\citenamefont{White}(1989)\citenamefont{White}}]{White1989}
\bibinfo{author}{\bibnamefont{White}, \bibfnamefont{E. L.}},
(\bibinfo{year}{1989}).
\emph{\bibinfo{title}{Cortical Circuits. Synaptic Organization of the Cerebral Cortex.}}
\bibinfo{Book}{Birkhauser, Boston 1989}.

\bibitem[{\citenamefont{Sengupta  et~al.}(2010)\citenamefont{Sengupta, Stemmler, Laughlin, and Niven}}]{Sengupta2010}
\bibinfo{author}{\bibnamefont{Sengupta}, \bibfnamefont{B.}},
\bibinfo{author}{\bibnamefont{Stemmler}, \bibfnamefont{M.}},
\bibinfo{author}{\bibnamefont{Laughlin}, \bibfnamefont{S.~B.}}, \bibnamefont{and}
\bibinfo{author}{\bibnamefont{Niven}, \bibfnamefont{J.~E.}}, (\bibinfo{year}{2010}).
\emph{\bibinfo{title}{Action potential energy efficiency varies among neuron types in vertebrates and invertebrates.}}
\bibinfo{journal}{PLoS Computational Biology} \textbf{\bibinfo{volume}{6}},
\bibinfo{pages}{e1000840}.


\bibitem[{\citenamefont{Carter and Bean}(2009)}]{Carter2009}
\bibinfo{author}{\bibnamefont{Carter}, \bibfnamefont{B.~C.}} \bibnamefont{and}
\bibinfo{author}{\bibnamefont{Bean}, \bibfnamefont{B.~P.}}, (\bibinfo{year}{2009}).
\emph{\bibinfo{title}{Sodium Entry during Action Potentials of Mammalian Neurons: Incomplete Inactivation and Reduced Metabolic Efficiency in Fast-Spiking Neurons.}}
\bibinfo{journal}{Neuron} \textbf{\bibinfo{volume}{64}}, \bibinfo{pages}{898-909}.

\bibitem[{\citenamefont{Crotty  et~al.}(2006)\citenamefont{Crotty, Sangrey, and Levy}}]{Crotty2006}
\bibinfo{author}{\bibnamefont{Crotty}, \bibfnamefont{P.}},
\bibinfo{author}{\bibnamefont{Sangrey}, \bibfnamefont{T.}}, \bibnamefont{and}
\bibinfo{author}{\bibnamefont{Levy}, \bibfnamefont{W.~B.}}, (\bibinfo{year}{2006}).
\emph{\bibinfo{title}{Metabolic energy cost of action potential velocity.}}
\bibinfo{journal}{J. Neurophysiol} \textbf{\bibinfo{volume}{96}},
\bibinfo{pages}{1237-1246}.


\bibitem[{\citenamefont{Ames III}(2000)}]{Ames2000}
\bibinfo{author}{\bibnamefont{Ames III}, \bibfnamefont{A.}}, (\bibinfo{year}{2000}).
\emph{\bibinfo{title}{CNS energy metabolism as related to function.}} \bibinfo{journal}{Brain Research Reviews} \textbf{\bibinfo{volume}{34}}, \bibinfo{pages}{42-68}.



\bibitem[{\citenamefont{Kandel et~al.}(1991)\citenamefont{Kandel, Schwartz, and
Jessell}}]{Kandel1991}
\bibinfo{author}{\bibnamefont{Kandel}, \bibfnamefont{E. R.}},
\bibinfo{author}{\bibnamefont{Schwartz}, \bibfnamefont{J. H.}}, \bibnamefont{and}
\bibinfo{author}{\bibnamefont{Jessell}, \bibfnamefont{T. M.}}, (\bibinfo{year}{1991}).
\emph{\bibinfo{title}{Principles of neural science (3rd ed.).}}
\bibinfo{Book}{Norwalk: Appleton and Lange}.

\bibitem[{\citenamefont{Laughlin et~al}(1998)\citenamefont{Laughlin,
   Steveninck and Anderson}}]{Laughlin1998}
\bibinfo{author}{\bibnamefont{Laughlin}, \bibfnamefont{S.~B.}},
\bibinfo{author}{\bibnamefont{Steveninck}, \bibfnamefont{R.~R. de Ruyter van}}, \bibnamefont{and}
\bibinfo{author}{\bibnamefont{Anderson}, \bibfnamefont{J.C.}}, (\bibinfo{year}{1998}).
\emph{\bibinfo{title}{The metabolic cost of neural information.}}
\bibinfo{journal}{Nature Neuroscience} \textbf{\bibinfo{volume}{1}}, \bibinfo{pages}{36-41}.

\bibitem[{\citenamefont{Alle et~al.}(2009)\citenamefont{Alle, Roth, and
Geiger}}]{Alle2009}
\bibinfo{author}{\bibnamefont{Alle}, \bibfnamefont{H.}},
\bibinfo{author}{\bibnamefont{Roth}, \bibfnamefont{A.}}, \bibnamefont{and}
\bibinfo{author}{\bibnamefont{Geiger}, \bibfnamefont{J.~R.~P.}}, (\bibinfo{year}{2009}).
\emph{\bibinfo{title}{Energy-efficient action potentials in hippocampal mossy fibers.}}
\bibinfo{journal}{Science} \textbf{\bibinfo{volume}{325}},
\bibinfo{pages}{1405-1408}.


\bibitem[{\citenamefont{Attwell and Laughlin}(2001)}]{Attwell2001}
\bibinfo{author}{\bibnamefont{Attwell}, \bibfnamefont{D.}} \bibnamefont{and}
\bibinfo{author}{\bibnamefont{Laughlin}, \bibfnamefont{S.~B.}}, (\bibinfo{year}{2001}).
\emph{\bibinfo{title}{An energy budget for signaling in the grey matter of the brain.}}
\bibinfo{journal}{Journal of Cerebral Blood Flow and Metabolism}
\textbf{\bibinfo{volume}{21}}, \bibinfo{pages}{1133-1145}.

\bibitem[{\citenamefont{Hertz et~al.}(2013)\citenamefont{Hertz, Junnan, Dan, Enzhi, Li and Liang}}]{Hertz2013}
\bibinfo{author}{\bibnamefont{Hertz}, \bibfnamefont{L.}},
\bibinfo{author}{\bibnamefont{Junnan}, \bibfnamefont{X.}},
\bibinfo{author}{\bibnamefont{Dan}, \bibfnamefont{S.}},
\bibinfo{author}{\bibnamefont{Enzhi}, \bibfnamefont{Y.}},
\bibinfo{author}{\bibnamefont{Li}, \bibfnamefont{G.}}, \bibnamefont{and}
\bibinfo{author}{\bibnamefont{Liang}, \bibfnamefont{P.}}, (\bibinfo{year}{2013}).
\emph{\bibinfo{title}{Astrocytic and neuronal accumulation of elevated extracellular K+ with a 2/3 K+/Na+ flux ratio - consequences for energy metabolism, osmolarity and higher brain function}}
\bibinfo{journal}{Front. Comput. Neurosci. } \textbf{\bibinfo{volume}{7}},
\bibinfo{number}{114}, \bibinfo{doi}{doi:10.3389/fncom.2013.00114}.
%


\bibitem[{\citenamefont{Moujahid et~al}(2011)\citenamefont{Moujahid,
d'Anjou, Torrealdea and Torrealdea}}]{Moujahid2011}
\bibinfo{author}{\bibnamefont{Moujahid}, \bibfnamefont{A.}},
\bibinfo{author}{\bibnamefont{d'Anjou}, \bibfnamefont{A.}},
\bibinfo{author}{\bibnamefont{Torrealdea}, \bibfnamefont{F.~J.}}, \bibnamefont{and}
\bibinfo{author}{\bibnamefont{Torrealdea}, \bibfnamefont{F.}}, (\bibinfo{year}{2011}).
\emph{\bibinfo{title}{Energy and information in Hodgkin-Huxley neurons.}}
\bibinfo{journal}{Physical Review E}
\textbf{\bibinfo{volume}{83}}, \bibinfo{pages}{031912}.


\bibitem[{\citenamefont{Moujahid et~al}(2012)\citenamefont{Moujahid and d'Anjou}}]{Moujahid2012}
\bibinfo{author}{\bibnamefont{Moujahid}, \bibfnamefont{A.}}, \bibnamefont{and}
\bibinfo{author}{\bibnamefont{d'Anjou}, \bibfnamefont{A.}}, (\bibinfo{year}{2012}).
\emph{\bibinfo{title}{Metabolic efficiency with fast spiking in the squid axon.}}
\bibinfo{journal}{Front. Comput. Neurosci.}
\textbf{\bibinfo{volume}{6}}, \bibinfo{number}{95}, \bibinfo{DOI}{doi: 10.3389/fncom.2012.00095}.

\bibitem[{\citenamefont{Hodgkin}(1975)}]{Hodgkin1975}
\bibinfo{author}{\bibnamefont{Hodgkin}, \bibfnamefont{A.~L.} },
  \bibinfo{journal}{Philosophical Transactions of the Royal Society B:
  Biological Sciences} \textbf{\bibinfo{volume}{270}}, \bibinfo{pages}{297}
  (\bibinfo{year}{1975}).

\bibitem[{\citenamefont{Lennie}(2003)}]{Lennie2003}
\bibinfo{author}{\bibnamefont{Lennie}, \bibfnamefont{P.}}, (\bibinfo{year}{2003}).
\emph{\bibinfo{title}{The Cost of Cortical Computation.}}
\bibinfo{journal}{Current Biology} \textbf{\bibinfo{volume}{13}},
\bibinfo{pages}{493-497}.



\bibitem[{\citenamefont{Moujahid et~al.}(2010)\citenamefont{Moujahid, D'Anjou, F. Torrealdea
  and F.J. Torrealdea}}]{Moujahid2010}
\bibinfo{author}{\bibnamefont{Moujahid}, \bibfnamefont{A.}},
  \bibinfo{author}{\bibnamefont{D'Anjou, A.}},
  \bibinfo{author}{\bibnamefont{Torrealdea, F.}} \bibnamefont{and}
  \bibinfo{author}{\bibnamefont{Torrealdea}, \bibfnamefont{F.J.}},
  \bibinfo{journal}{Sixth International Conference on Natural Computation (ICNC)} \textbf{\bibinfo{volume}{6}}, \bibinfo{pages}{3156-3163}
  (\bibinfo{year}{2010}).


\bibitem[{\citenamefont{Stein}(2002)}]{Stein2002}
\bibinfo{author}{\bibnamefont{Stein}, \bibfnamefont{W.~D.}},
  \bibinfo{journal}{International Review of Cytology}
  \textbf{\bibinfo{volume}{215}}, \bibinfo{pages}{231} (\bibinfo{year}{2002}).


\bibitem[{\citenamefont{Buzsaki et~al.}(2007)\citenamefont{Buzsaki, K.Kaila,
  and Raichle}}]{Buzsaki}
\bibinfo{author}{\bibnamefont{Buzsaki}, \bibfnamefont{G.}},
  \bibinfo{author}{\bibnamefont{Kaila, K.}}, \bibnamefont{and}
  \bibinfo{author}{\bibnamefont{Raichle}, \bibfnamefont{M.}},
  \bibinfo{journal}{Neuron} \textbf{\bibinfo{volume}{56}}, \bibinfo{pages}{771}
  (\bibinfo{year}{2007}).

\bibitem[{\citenamefont{Pospischil et~al.}(2008)\citenamefont{Pospischil, Toledo-Rodriguez,
Monier, Piwkowska, Bal, Frégnac, Markram and Destexhe}}]{Pospischil2008}
\bibinfo{author}{\bibnamefont{Pospischil}, \bibfnamefont{M.}},
\bibinfo{author}{\bibnamefont{Toledo-Rodriguez}, \bibfnamefont{M.}},
\bibinfo{author}{\bibnamefont{Monier}, \bibfnamefont{C.}},
\bibinfo{author}{\bibnamefont{Piwkowska}, \bibfnamefont{Z.}},
\bibinfo{author}{\bibnamefont{Bal}, \bibfnamefont{T.}},
\bibinfo{author}{\bibnamefont{Frégnac}, \bibfnamefont{Y.}},
\bibinfo{author}{\bibnamefont{Markram}, \bibfnamefont{H.}}, \bibnamefont{and}
\bibinfo{author}{\bibnamefont{Destexhe}, \bibfnamefont{A.}}, (\bibinfo{year}{2008}).
\emph{\bibinfo{title}{Minimal Hodgkin-Huxley type models for different classes of cortical and thalamic neurons.}}
\bibinfo{journal}{Biological Cybernetics} \textbf{\bibinfo{volume}{99}},
\bibinfo{pages}{427-441}.

\bibitem[{\citenamefont{Guo et~al}(2008)\citenamefont{Guo,
Rubin, McIntyre, Vitek and Terman}}]{Guo2008}
\bibinfo{author}{\bibnamefont{Guo}, \bibfnamefont{Y.}},
\bibinfo{author}{\bibnamefont{Rubin}, \bibfnamefont{J. E.}},
\bibinfo{author}{\bibnamefont{McIntyre}, \bibfnamefont{C. C.}},
\bibinfo{author}{\bibnamefont{Vitek}, \bibfnamefont{J.L.}}, \bibnamefont{and}
\bibinfo{author}{\bibnamefont{Terman}, \bibfnamefont{D.}}, (\bibinfo{year}{2008}).
\emph{\bibinfo{title}{Thalamocortical Relay Fidelity Varies Across Subthalamic Nucleus Deep Brain Stimulation Protocols in a Data-Driven Computational Model.}}
\bibinfo{journal}{Journal of Neurophysiology}
\textbf{\bibinfo{volume}{99}}, \bibinfo{number}{3}, \bibinfo{pages}{1477-1492}.


\bibitem[{\citenamefont{Gibbons and Chakraborti}(1996)\citenamefont{Gibbons and Chakraborti}}]{Gibbons2011}
\bibinfo{author}{\bibnamefont{Gibbons}, \bibfnamefont{J. D.}}, \bibnamefont{and}
\bibinfo{author}{\bibnamefont{Chakraborti}, \bibfnamefont{S.}}, (\bibinfo{year}{2011}).
\emph{\bibinfo{title}{Nonparametric Statistical Inference.}}
\bibinfo{Book}{5th Ed., Boca Raton, FL: Chapman and Hall-CRC Press, Taylor and Francis Group}


\bibitem[{\citenamefont{Wang and Buzsáki}(1996)\citenamefont{Wang and Buzsáki}}]{Wang1996}
\bibinfo{author}{\bibnamefont{Wang}, \bibfnamefont{X-J.}}, \bibnamefont{and}
\bibinfo{author}{\bibnamefont{Buzsáki}, \bibfnamefont{G.}}, (\bibinfo{year}{1996}).
\emph{\bibinfo{title}{Gamma Oscillation by Synaptic Inhibition in a Hippocampal Interneuronal Network Model.}}
\bibinfo{journal}{Journal of Neuroscience}
\textbf{\bibinfo{volume}{16}}, \bibinfo{number}{20}, \bibinfo{pages}{6402-6413}.

\bibitem[{\citenamefont{Jansen et~al}(2003)\citenamefont{Jansen,
Shen, Zhang, Wolkowicz and Balschi}}]{Jansen2003}
\bibinfo{author}{\bibnamefont{Jansen}, \bibfnamefont{Maurits A.}},
\bibinfo{author}{\bibnamefont{Shen}, \bibfnamefont{H.}},
\bibinfo{author}{\bibnamefont{Zhang}, \bibfnamefont{L.}},
\bibinfo{author}{\bibnamefont{Wolkowicz}, \bibfnamefont{Paul E.}}, \bibnamefont{and}
\bibinfo{author}{\bibnamefont{Balschi}, \bibfnamefont{James A.}}, (\bibinfo{year}{2003}).
\emph{\bibinfo{title}{Energy requirements for the Na+ gradient in the oxygenated isolated heart: effect of changing the free energy of ATP hydrolysis.}}
\bibinfo{journal}{American Journal of Physiology - Heart and Circulatory Physiology}
\textbf{\bibinfo{volume}{285}}, \bibinfo{pages}{H2437--H2445}.


\bibitem[{\citenamefont{Ereciriska et~al}(1989)\citenamefont{Ereciriska and Silver}}]{Ereciriska1989}
\bibinfo{author}{\bibnamefont{Ereciriska}, \bibfnamefont{M.}}, \bibnamefont{and}
\bibinfo{author}{\bibnamefont{Silver}, \bibfnamefont{Ian A.}}, (\bibinfo{year}{1989}).
\emph{\bibinfo{title}{ATP and Brain Function.}}
\bibinfo{journal}{Journal of Cerebral Blood Flow and Metabolism}
\textbf{\bibinfo{volume}{9}}, \bibinfo{pages}{2-19}.

%
\bibitem[{\citenamefont{Hodgkin and Huxley}(1952)}]{Hodgkin1952}
\bibinfo{author}{\bibnamefont{Hodgkin}, \bibfnamefont{A.~L.}} \bibnamefont{and}
\bibinfo{author}{\bibnamefont{Huxley}, \bibfnamefont{A.~F.}}, (\bibinfo{year}{1952}).
\emph{\bibinfo{title}{A quantitative description of membrane current and its application to conduction and excitation in nerve.}}
\bibinfo{journal}{The Journal of Physiology} \textbf{\bibinfo{volume}{117}}, \bibinfo{pages}{500-544}.


\bibitem[{\citenamefont{Connors}(1990)}]{Connors1990}
\bibinfo{author}{\bibnamefont{Connors}, \bibfnamefont{Barry W.}}, \bibnamefont{and}
\bibinfo{author}{\bibnamefont{Gutnick}, \bibfnamefont{Michael J.}}, (\bibinfo{year}{1990}).
\emph{\bibinfo{title}{Intrinsic firing patterns of diverse neocortical neurons.}}
\bibinfo{journal}{Trends in Neurosciences}
\textbf{\bibinfo{volume}{13}}, \bibinfo{pages}{99-104}.



\bibitem[{\citenamefont{Koch}(1998)}]{Koch1998}
\bibinfo{author}{\bibnamefont{Koch}, \bibfnamefont{C.}}, (\bibinfo{year}{1998}).
\emph{\bibinfo{title}{Biophysics of computation: Information Processing in Single Neurons.}} \bibinfo{Book}{New York: Oxford University Press}.

\bibitem[{\citenamefont{Falk}(1990)}]{Falk1990}
\bibinfo{author}{\bibnamefont{Falk}, \bibfnamefont{D.}}, (\bibinfo{year}{1990}).
\emph{\bibinfo{title}{Brain evolution in Homo: the radiator theory.}} \bibinfo{journal}{Behavioral and Brain Sciences} \textbf{\bibinfo{volume}{13}}, \bibinfo{pages}{333-381}.

\bibitem[{\citenamefont{Kiyatkin}(2007)}]{Kiyatkin2007}
\bibinfo{author}{\bibnamefont{Kiyatkin}, \bibfnamefont{E.~A.}}, (\bibinfo{year}{2007}).
\emph{\bibinfo{title}{Brain temperature fluctuations during physiological and pathological conditions.}} \bibinfo{journal}{European Journal of Applied Physiology} \textbf{\bibinfo{volume}{101}}, \bibinfo{pages}{3-17}.




\bibitem[{\citenamefont{Sohal et al.}(2002)}]{Sohal2002}
\bibinfo{author}{\bibnamefont{Sohal}, \bibfnamefont{Vikaas S.}}, \bibnamefont{and}
\bibinfo{author}{\bibnamefont{Huguenard}, \bibfnamefont{John R.}}, (\bibinfo{year}{2002}).
\emph{\bibinfo{title}{Reciprocal inhibition controls the oscillatory state in thalamic networks.}}
\bibinfo{journal}{Neurocomputing}
\textbf{\bibinfo{volume}{44-46}}, \bibinfo{pages}{653-659}.









\end{thebibliography}
\end{document}